\documentclass{article}
\usepackage{amsmath}
\usepackage{graphicx}
\usepackage{amsfonts}
\title{How to Compute a Puiseux Expansion}
\author{Nicholas J. Willis, Annie K. Didier, Kevin M. Sonnanburg \\ Whitworth University}
\begin{document}
\maketitle

\begin{abstract} 
In this paper, an explanation of the Newton-Peiseux algorithm is given. This
explanation is supplemented with well-worked and explained examples of how to
use the algorithm to find fractional power series expansions for all branches of a
polynomial at the origin.
\end{abstract}
\section{Introduction}
Given a polynomial $f(x, y)$, you can always solve for $y$ in terms of $x$ by means of
a fractional power series
\begin{equation*}
y=c_1x^{\gamma_1}+c_2x^{\gamma_1+\gamma_2}+c_3x^{\gamma_1+\gamma_2+\gamma_3}+\cdots
\end{equation*}
The method for doing this is often explained in books in the context of a proof of the
algebraic closure of the field of fractional power series (e.g., [1]).
Our purpose here is not to discuss the proof of such a theorem, but to explain
clearly, with several examples, the algorithm for obtaining these expansions. \\
Here is a general description of the algorithm.\\
\begin{description}
\item[1.] Given $f(x, y) = 0$ draw the Newton polygon of $f(x, y)$.  This is done by plotting a point on
$\mathbb{R}^2$ for each term of $f$ with the term $kx^ay^b$ being plotted to the point $(b,a)$  The Newton polygon is the smallest convex shape that contains all the points plotted.
\item[2.] Take a segment of the Newton polygon from the set of segments where each point plotted is on, above, or to the right of the segments.  
\item[3.] The first exponent $\gamma_1$ will be the negative of the slope of this segment.
\item[4.] Find $f(x,x^{\gamma_1}(c_1+ y_1))$.
\item[5.] Take the lowest terms in $x$ alone. Since $f(x,y) = 0$, these must cancel and so
you can solve for $c_1$,
\item[6.] Taking the values of $\gamma_1$ and $c_1$ and $\beta =$ the '$x$-intercept' on the Newton polygon
of the segment we have chosen, we now find $f_1(x, y_1)$.
\item[7.] $f_1(x,y_1) = x^{-\beta}f(x,x^{\gamma_1}(c_1+y_1))$
\item[8.] Now repeat the process for $f_1(x,y_1)$ to find $\gamma_2$ and $c_2$.
\item[9.] Continue this process until one of two things happen \end{description} 
\begin{itemize}
\item $f_n(x, y_n)$ has a factor of $y_n$.
\item The Newton Polygon of $f_n(x, y_n)$ consists of a single segment with only two
vertices, one on each axis. 
\end{itemize}

If the former occurs, either the $y$ factor(s) can be factored out and you can continue,
or (if you are not left with a segment in the Newton polygon) the series terminates.
If the latter occurs, plug in the rest of the series, solving for coefficients by letting lowest
powers of $x$ cancel. Use the denominator of the slope of the segment to determine the
increase in each $\gamma_i$, (e.g., if the slope is $\frac{3}{2}$ , each $\gamma_i$, should be $\frac{1}{2}$) . This will be greatly
clarified by the following examples.

Any polynomial of the form $f(x, y) = 0$ can be solved for $y$ in terms of $x$ explicitly
around a point. Since by simple translation any point on a curve can be moved
to the origin, we will only expand a solution of $f(x,y) = 0$ at the origin. The
Newton-Puiseux algorithm gives us that any solution of $f(x, y) = 0$ has the form
$y=c_1x^{\gamma_1}+c_2x^{\gamma_1+\gamma_2}+c_3x^{\gamma_1+\gamma_2+\gamma_3}+\cdots$, where $\gamma_i \in \mathbb{Q}^+$ and the $c_i \in \mathbb{Q}^+$. Any time $\gamma_i$
is referred to, we are speaking of an exponent, and anytime a $c_i$ is referred to, we are
referring to a coefficient.
\section{Examples}
\emph{Example 1:}
In this example we show a more simple way to compute the expansion that can be used if the second case in step 9, mentioned above, occurs where the Newton Polygon has a horizontal intercept of one.
\begin{figure} [ht]
\centering
\caption{Newton Polygon for $f(x,y)$}
\scalebox{.7}{\includegraphics{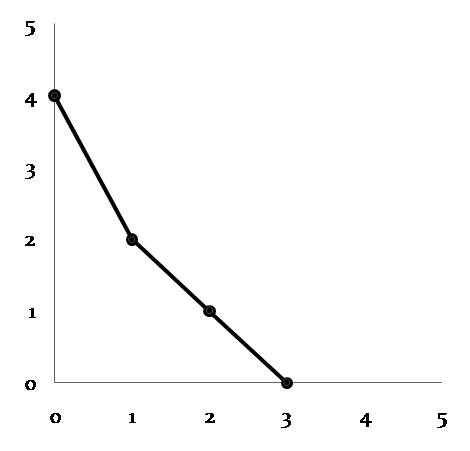}}
\end{figure}
\medskip
\begin{equation}
f(x,y)=2x^4+x^2y+4xy^2+4y^3=0 \label{1}
\end{equation}
\\
We know that this polynomial has solutions of the form 
\begin{equation}
y=c_1x^{\gamma_1}+c_2x^{\gamma_1+\gamma_2}+c_3x^{\gamma_1+\gamma_2+\gamma_3}+\cdots. \label{2}
\end{equation}
The possible choices for $\gamma_1$ correspond to the negatives of the slopes of the lower segments of the Newton Polygon (Figure 1). The slopes of the lower segments of the Newton Polygon corresponding to (1) are -2, and -1. We will start by considering the case where $\gamma_1=-(-2)=2$.

Let $\gamma_1=2$. Then by factoring $x^{\gamma_1}$ (or $x^2$) out of (2) we obtain, 
\begin{equation}
y=x^2(c_1+y_1) \label{3}, 
\end{equation}
where $c_1$ is the first coefficient in the series and $y_1$ is
\begin{equation}
y_1=c_2x^{\gamma_2}+c_3x^{\gamma_2+\gamma_3}+c_4x^{\gamma_2+\gamma_3+\gamma_4}+\cdots \label{4}
\end{equation}

i.e. the rest of the series.

So, substituting (3) into (1) we get, 
\begin{equation}
2x^4+x^4(c_1+y_1)+4x^5(c_1+y_1)^2+4x^6(c_1+y_1)^3=0 \label{5}
\end{equation}

The vertical intercept of the segment we used to find $f_1$ on the Newton Polygon, which we call $\beta$, gives the terms of lowest degree in $x$ alone. These terms must sum to zero because $f(x,y)=0$. In this case, $\beta=4$, so \\ $2x^4+c_1x^4=0$ and solving for $c_1$ we obtain $c_1=-2$. 
Substituting  $c_1$ into (5) and dividing by $x^4$ gives, 

\begin{center}
$x^{-4}f(x,x^2(c_1+y_1))$ $=2+(-2+y_1)+4x(-2+y_1)^2+4x^2(-2+y_1)^3$
\end{center}
\begin{equation}
                                                =y_1+16x-16xy_1+4xy_1^2-32x^2+48x^2y_1-24x^2y_1^2+4x^2y_1^3 \label{6}
 \end{equation}
\pagebreak
\begin{figure} [ht]
\centering
\caption{Newton Polygon for $f_1(x,y_1)$}
\scalebox{.7}{\includegraphics{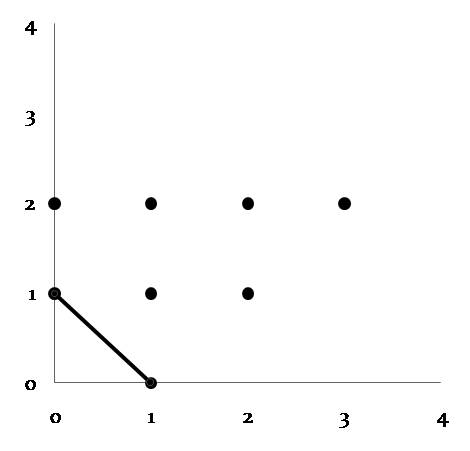}}
\end{figure}

We call this new polynomial $f_1(x,y_1)$, so now that we have our new function $f_1(x,y_1)$ our next goal is to find $c_2$. We find the new Newton Polygon for $f_1(x,y_1)$ (Figure 2), from which we get $\beta_2=1$ and $\gamma_2=1$. Again, we let the terms of lowest degree in $x$ alone equal zero. We have now reached a point in our algorithm where the calculations become much
simpler. Notice that in the previous Newton polygon the intercept with the horizontal axis was one and we had only one segment. This
tells us that we can calculate the rest of the $c_i$ and $\gamma_i$ directly from the current function, in this case $f_1$.This is true because when we obtain such a Newton Polygon, we know that the powers of $x$ are going to increase by at least 1/(denominator of $\gamma_2$) from now on (with the possibility that $c_i$ could be = 0 ); in other words,
$\gamma_3 = \gamma_4 = \cdots= \gamma_n = \cdots =$.
Since $\gamma_2$, we know that
\begin{center}
$y_1=c_2x+c_3x^2+c_4x^3+c_5x^4+\cdots$
\end{center}
So, substituting $y_1$ into (6) and letting the lowest degree terms of x cancel we see 
\begin{center}
$c_2x+16x=0$  

so, $c_2=-16$
\end{center}
We find $c_3$ in the same manner. To do this, we let the $x^2$ terms equal zero. Substituting $y_1$ into (6) we see that,
\begin{center}
$c_3x^2-16c_2x^2-32x^2=0$

$c_3+256-32=0$ 

so, $c_3=-224$
\end{center}
Now, we can let the $x^3$ terms in $f_1$ be equal to zero. So, 
\begin{center}
$c_4x^3-16c_3x^3+4c_2^2x^3+48c_2x^3=0$

$c_4+3584+1024-768=0$ 

so, $c_4=-3840$
\end{center}
Next, we let the $x^4$ terms equal zero.
\begin{center}
$c_5x^4-16c_4x^4+8c_2c_3x^4+48c_3x^4-24c_1^2x^4=0$

$c_5+61440+28672-11184-10752=0$ 

so, $c_5=-39504$
\end{center}
This process can be continued to calculate more terms, but we will stop with $c_5$. So, one explicit solution of the polynomial is
\begin{center}
$y=-2x^2-16x^3-224x^4-3840x^5-39504x^6+\cdots$
\end{center}
\bigskip
We have examined the case where $\gamma_1=2$, and now we must return to the original polynomial $f(x,y)$ and examine the case where $\gamma_1=1$. Recall,\\ $f(x,y)=2x^4+x^2y+4xy^2+4y^3=0$. Let $\gamma_1=1$. Then, $y=x(c_1+y_1)$ and substituting $\gamma_1$ and $y$ into (1) we obtain 
\begin{equation}
2x^4+x^3(c_1+y_1)+4x^3(c_1+y_1)^2+4x^3(c_1+y_1)^3=0 \label{7}
\end{equation}
As in the case where $\gamma_1=2$, 
\begin{center}
$y_1=c_2x^{\gamma_2}+c_3x^{\gamma_2+\gamma_3}+c_4x^{\gamma_2+\gamma_3+\gamma_4}+\cdots$  
\end{center}
From the Newton Polygon (Figure 1) we see that $\beta_1=3$. Therefore, we can let the $x^3$ terms equal zero.
\begin{center}
$c_1x^3+4c_1^2x^3+4c_1^2x^3=0$  
\\
$c_1(1+4c_1+4c_1^2)=0$
\\
$c_1(1+2c_1)^2=0$
\\
$c_1=0, -\frac{1}{2}, -\frac{1}{2}$ 
\\
$c_1=0$ is trivial, so assume, $c_1=-\frac{1}{2}$.
\end{center}
Substituting $c_1$ into (7) and dividing by $x^3$ gives: 
\begin{align*}
x^{-3}f(x,x(c_1+y_1)) &=2x+(-\frac{1}{2}+y_1)+4(-\frac{1}{2}+y_1)^2+4(-\frac{1}{2}+y_1)^3\\
                              &=2x+y_1-4y_1+4y_1^2++3y_1-6y_1^2+4y_1^3\\
f_1(x,y)                   &=2x-2y_1^2+4y_1^3\\
\end{align*}
$\gamma_2=2$, so
\begin{center}
$y_1=x^\frac{1}{2}(c_2+y_2)$ so, $2x-2x(c_2+y_2)^2+4x^\frac{3}{2}(c_2+y_2)^3=0$
\end{center}
where $c_2$ is the second coefficient and $y_2$ is the rest of the series. Substituting $y_1$ into $f_1(x,y)$ yields
\begin{center}
$2x-2x(c_2+y_2)^2+4x^\frac{3}{2}(c_2+y_2)^3=0$
\end{center}

From the new Newton Polygon (Figure 1) we see that $\beta_2=1$, so we can let the $x$ terms equal zero.
\begin{center}
 $2x-2xc_2^2=0$ 
\end{center}
\begin{equation}
 c_2=\pm1 \label{8}
\end{equation}
Let $c_2=1$. (We must return later and consider the case where $c_2=-1$.) Then substituting $c_2$ into $f_1(x,y)$ and dividing by $x$ we get
\begin{align*}
x^{-1}f(x,y)_1 &=2-2(1+y_2)^2+4x^\frac{1}{2}(1+y_2)^3
\\
&=-4y_2-2y_2^2+4x^\frac{1}{2}+12x^\frac{1}{2}y_2+12x^\frac{1}{2}y_2^2+4x^{1}{2}y_2^3
\end{align*}
Where $y_2=c_3x^{\gamma_3}+c_4x^{\gamma_3+\gamma_4}+c_5x^{\gamma_3+\gamma_4+\gamma_5}+\cdots$
\\
The horizontal intercept is one, so $y_2=c_3x^\frac{1}{2}+c_4x+c_5x^\frac{3}{2}+\cdots$
\\
Letting the lowest degree terms of $x$ ($x^\frac{1}{2}$) equal zero we find
\begin{center}
$-4c_3+4=0$ 
\\
$c_3=1$ 
\end{center}
Now,  we can let the $x$ terms equal zero.
\begin{center}
$-4c_4x-2c_3^2x+12c_3x=0$
\\
$-4c_4-2+12=0$ 
\\
so, $c_4=\frac{5}{2}$
\end{center}
We now let the $x^\frac{3}{2}$ terms equal zero.
\begin{center}
$-4c_5x^\frac{3}{2}-4c_3c_4x^\frac{3}{2}+12c_4\frac{3}{2}+12c_3^2x^\frac{3}{2}=0$

$-4c_5-10+30+12=0$ 

so, $c_5=8$
\end{center}
As in the case where $\gamma_1=2$ we could continue this process and calculate more terms. Stopping here, another explicit solution of $f(x,y)$ is
\begin{center}
$y=-\frac{1}{2}x+x^\frac{3}{2}+x^2+\frac{5}{4}x^\frac{5}{2}+\frac{27}{4}x^3+\cdots$
\end{center}
\bigskip

Now, we must return to (8) and consider the case where $c_1=-1$. Substituting $c_1$ into $f_1(x,y)$ and dividing by $x^4$ gives
\begin{align*}
x^{-4} f_1(x,y) &=2-2(-1+y_2)^2+4x^\frac{1}{2}(-1+y_2)^3
\\
&=4y_2-2y_2^2-4x^\frac{1}{2}+12x^\frac{1}{2}y_2-12x^\frac{1}{2}y_2^2+4x^\frac{1}{2}y_2^3
\end{align*}
Let $y_2=c_3x^{\gamma_3}+c_4x^{\gamma_3+\gamma_4}+c_5x^{\gamma_3+\gamma_4+\gamma_5}+\cdots$.
$\gamma_3=\frac{1}{2}$, so
$y_2=c_3x^\frac{1}{2}+c_4x+c_x^\frac{3}{2}+\cdots$.
Letting the $x^\frac{1}{2}$ terms equal zero, 
\begin{center}
$4c_3x^\frac{1}{2}+4x^\frac{1}{2}=0$
\\
so, $c_3=-1$
\end{center}
Letting the $x$ terms equal zero,
\begin{center}
$4c_4x-c_3^2x+12c_3x=0$
\\
$4c_4-1-12=0$ 
\\
so, $c_4=\frac{13}{4}$
\end{center}
Now we let the $x^\frac{3}{2}$ terms equal zero and solve for $c_5$,
\begin{center}
$4c_5+12c_4-12c_3^2=0$
\\
$4c_5+39-12=0$ 
\\
so, $c_5=-\frac{27}{4}$
\end{center}
So, we find that the last explicit solution for $f(x,y)=0$ is
\begin{center}
 $y=-\frac{1}{2}x-x^\frac{3}{2}-x^2+\frac{13}{4}x^\frac{5}{2}-\frac{27}{4}x^3$
\end{center}
So, the three solutions for $f(x,y)=0$ are
\begin{center}
$y=-2x^2-16x^3-224x^4-3840x^5-39504x^6+\cdots$

$y=-\frac{1}{2}x+x^\frac{3}{2}+x^2+\frac{5}{4}x^\frac{5}{2}+\frac{27}{4}x^3+\cdots$

$y=-\frac{1}{2}x-x^\frac{3}{2}-x^2+\frac{13}{4}x^\frac{5}{2}-\frac{27}{4}x^3$
\end{center}
\bigskip
\emph{Example 2:}
In this example we see that the Puiseux algorithm can be used to compute highly accurate fractional powers and expansions with imaginary and complex terms.
\begin{figure} [ht]
\centering
\caption{Newton Polygon for $f(x,y)$}
\scalebox{.7}{\includegraphics{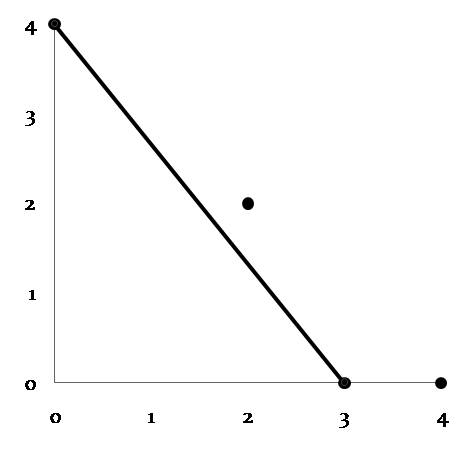}}
\end{figure}
\medskip
\begin{equation}
f(x,y)=x^5+8x^4-2x^2y^2-y^3+2y^4=0 \tag{1}
\end{equation}
We know that this polynomial has solutions of the form 
\begin{equation}
y=c_1x^{\gamma_1}+c_2x^{\gamma_1+\gamma_2}+c_3x^{\gamma_1+\gamma_2+\gamma_3}+\cdots. \tag{2}
\end{equation}
$\gamma_1$ corresponds to the negatives of the slope of the lower segment of the Newton Polygon (Figure 3). The slope of the lower segment of the Newton Polygon corresponding to (1) is $-\frac{4}{3}$. So, $\gamma_1=-(-\frac{4}{3})=\frac{4}{3}$.

Let $\gamma_1=\frac{4}{3}$ Then by factoring $x^{\gamma_1}$ (or $x^\frac{4}{3}$) out of (2) we obtain, 
\begin{equation}
y=x^\frac{4}{3}(c_1+y_1) \tag{3}, 
\end{equation}
where $c_1$ is the first coefficient in the series and
\begin{equation}
y_1=c_2x^{\gamma_2}+c_3x^{\gamma_2+\gamma_3}+c_4x^{\gamma_2+\gamma_3+\gamma_4}+\cdots \tag{4}
\end{equation}
\\
i.e. the rest of the series. \\
So, substituting (3) into (1) we get, 
\begin{equation}
x^5+8x^4-2x^\frac{14}{3}(c_1+y_1)^2-x^4(c_1+y_1)^3+2x^\frac{16}{3}(c_1+y_1)^4=0 \tag{5}
\end{equation}

The vertical intercept of the segment we use to find $\gamma_1$ on the Newton Polygon, which we call $\beta$, gives the terms of lowest degree in $x$ alone. These terms must sum to zero because $f(x,y)=0$. In this case, $\beta_1=4$ so, $8x^4-c_1^3x^4=0$ and solving for $c_1$ we obtain 
\begin{equation}
c_1=2, -1+i\sqrt{3}, -1-i\sqrt{3}   \tag{6}
\end{equation}
We assume $c_1=2$, but we must consider the other two cases later.
Substituting  $c_1$ into (5) and dividing by $x^4$ gives, 

\begin{align*}
x^{-4}f(x,x^2(c_1+y_1))&=f_1(x,y_1) \\
&=x-8x^\frac{2}{3}-8x^\frac{2}{3}y_1-2x^\frac{2}{3}y_1^2-12y_1-6y_1^2-y_1^3+32x^\frac{4}{3} \\
&+64x^\frac{4}{3}y_1+48x^\frac{4}{3}y_1^2+16x^\frac{4}{3}y_1^3+2x^\frac{4}{3}y_1^4 \tag{7}
\end{align*}

We call this new polynomial $f_1(x,y_1)$. So, now that we have our new function $f_1(x,y_1)$ our new goal is to find $c_2$. We find the new Newton Polygon (Figure 4), from which we get $\beta_2=\frac{2}{3}$ and $\gamma_2=\frac{2}{3}$. 
\begin{figure} [ht]
\centering
\caption{Newton Polygon for $f_1$}
\scalebox{.7}{\includegraphics{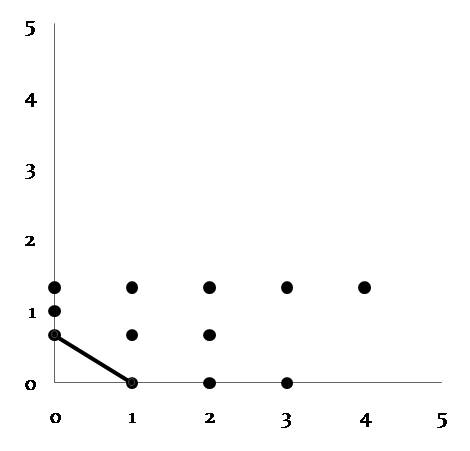}}
\end{figure}

Similar to (3), we factor $x^{\gamma_2}$ out of (2) and obtain
\begin{equation*}
y_1=x^\frac{2}{3}(c_2+y_2) \tag{8}
\end{equation*} 
where $c_2$ is the second coefficient in the series and 
\begin{equation*}
y_2=c_3x^{\gamma_3}+c_4x^{\gamma_3+\gamma_4}+c_5x^{\gamma_3+\gamma_4+\gamma_5}+\cdots
\end{equation*}
Subsituting (8) into $f_1$ we see that
\begin{align*}
0=&x-8x^{\frac{2}{3}}-8x^{\frac{4}{3}}(c_2+y_2)-2x^2(c_2+y_2)^2-12x^{\frac{2}{3}}(c_2+y_2)-6x^{4}{3}(c_2+y_2)^2 \\
&-x^2(c_2+y_2)^3+32x^{\frac{4}{3}}+64x^2(c_2+y_2)+48x^{\frac{8}{3}}(c_2+y_2)^2+16x^{\frac{10}{3}}(c_2+y_2)^3 \\
&+2x^4(c_2+y_2)^4 \tag{9} \\
\end{align*}
\pagebreak
\\
$\beta_2$ gives us the terms of lowest degree in $x$ alone, which we can again set equal to zero.
\begin{center}
$-8x^\frac{2}{3}-12c_2x^\frac{2}{3}=0$

so, $c_2=-\frac{2}{3}$
\end{center}
Substituting $c_2$ into (9) and dividing by $x^\frac{2}{3}$ we see
\begin{align*}
x^{-\frac{2}{3}}f_2(x,x^{\frac{2}{3}}(c_2+y_2))=&x^{\frac{1}{3}}+\frac{32}{81}x^{\frac{10}{3}}-\frac{1168}{27}x^{\frac{4}{3}}-\frac{128}{27}x^{\frac{8}{3}}+\frac{64}{3}x^2+\frac{104}{3}x^{\frac{2}{3}} \\
&+\frac{196}{3}x^{\frac{4}{3}}y_2-12y_2-6x^{\frac{2}{3}}y_2^2-x^{\frac{4}{3}}y_2^3-64x^2y_2 \\
&+48x^2y_2^2+\frac{64}{3}x^{\frac{8}{3}}y_2-32x^{\frac{8}{3}}y_2^2+16x^{\frac{8}{3}}y_2^3-\frac{64}{27}x^{\frac{10}{3}}y_2 \\
&+\frac{16}{3}x^{\frac{10}{3}}y_2^2-\frac{16}{3}x^{\frac{10}{3}}y_2^3+2x^{\frac{10}{3}}y_2^4 \tag{10} \\
\end{align*}
\begin{figure} [ht]
\centering
\caption{Newton Polygon for $f_2$}
\scalebox{.7}{\includegraphics{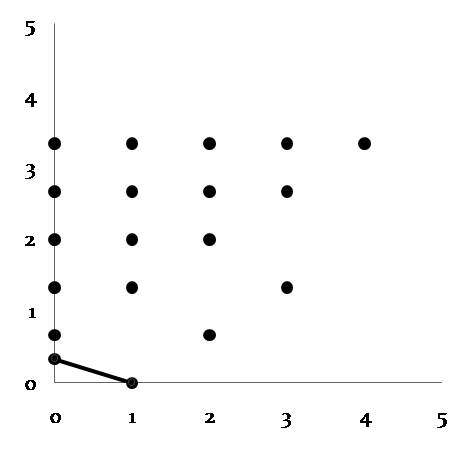}}
\end{figure}
\\
Now, we can see from Figure 5 that $\gamma_3=\frac{1}{3}$, so
\begin{equation}
y_2=x^\frac{1}{3}(c_3+y_3) \tag{11}
\end{equation}
where $c_3$ is the third coefficient in the series and 
\begin{equation*}
y_3=c_4x^{\gamma_4}+c_5x^{\gamma_4+\gamma_5}+c_6x^{\gamma_4+\gamma_5+\gamma_6}+\cdots 
\end{equation*}
We substitute $y_2$ into $f_2$ and find
\begin{align*}
0=&x^\frac{1}{3}+\frac{32}{81}x^\frac{10}{3}-\frac{1168}{27}x^\frac{4}{3}-\frac{128}{27}x^\frac{8}{3}+\frac{64}{3}x^2+\frac{104}{3}x^\frac{2}{3}+\frac{196}{3}x^\frac{5}{3} (c_3+y_3) \\ &-12x^\frac{1}{3}(c_3+y_3)-6x^\frac{4}{3}(c_3+y_3)^2-x^\frac{7}{3}(c_3+y_3)^3-64x^\frac{7}{3}(c_3+y_3) \\
&+48x^\frac{8}{3}(c_3+y_3)^2+\frac{64}{3}x^3(c_3+y_3)-32x^\frac{10}{3}(c_3+y_3)^2+16x^\frac{11}{3}(c_3+y_3)^3 \\
&-\frac{64}{27}x^\frac{11}{3}(c_3+y_3)+\frac{16}{3}x^4(c_3+y_3)^2-\frac{16}{3}x^\frac{13}{3}(c_3+y_3)^3+2x^\frac{14}{3}(c_3+y_3)^4 \\ \tag{12}
\end{align*}

Since $f_2=0$ the $x^\frac{1}{3}$ terms must also equal zero, so we can let
\begin{center}
$x^\frac{1}{3}-12c_3x^\frac{1}{3}=0$

so, $c_3=\frac{1}{12}$
\end{center}
Substituting $c_3$ into (12) and dividing by $x^\frac{1}{3}$ we see that
\begin{align*}
x^{-\frac{1}{3}}f_3&(x,x^\frac{1}{3}(c_3+y_3))=-\frac{9353}{216}x-\frac{1}{324}x^4-\frac{9217}{1728}x^2+\frac{14}{81}x^3+\frac{1}{27}x^\frac{11}{3}+\frac{16}{9}x^\frac{8}{3} \\
&-\frac{61}{324}x^\frac{10}{3}-\frac{119}{27}x^\frac{7}{3}+\frac{49}{9}x^\frac{4}{3}+\frac{104}{3}x^\frac{1}{3}-12y_3+\frac{64}{3}x^\frac{5}{3}+\frac{1}{10368}x^\frac{13}{3}\\ 
&-\frac{16}{3}x^3y_3-\frac{1}{9}x^4y_3-\frac{4}{3}x^4y_3^2+\frac{196}{3}x^\frac{4}{3}y_3+8x^\frac{7}{3}y_3+48x^\frac{7}{3}y_3^2+\frac{64}{3}x^\frac{8}{3}y_3 \\
&-\frac{55}{27}x^\frac{10}{3}y_3+4x^\frac{10}{3}y_3^2+\frac{8}{9}x^\frac{11}{3}y_3+\frac{16}{3}x^\frac{11}{3}y_3^2+\frac{1}{216}x^\frac{13}{3}y_3+\frac{1}{12}x^\frac{13}{3}y_3^2 \\
&+\frac{2}{3}x^\frac{13}{3}y_3^3-xy_3-6xy_3^2-\frac{3073}{48}x^2y_3-\frac{1}{4}x^2y_3^2-x^2y_3^3-32x^3y_3^2\\
&+16x^\frac{10}{3}y_3^3-\frac{16}{3}x^4y_3^3+2x^\frac{13}{3}y_3^4 \tag{13} \\
\end{align*}
\begin{figure} [ht]
\centering
\caption{Newton Polygon for $f_3$}
\scalebox{.7}{\includegraphics{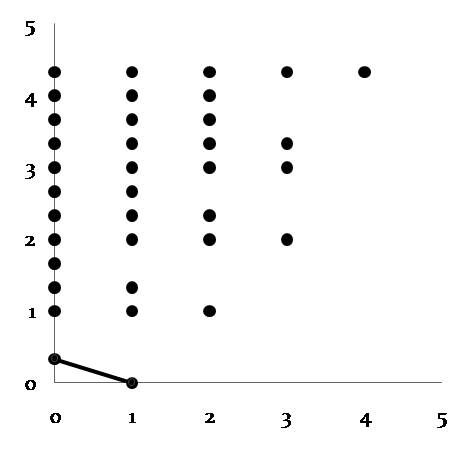}}
\end{figure}
\\
Now, we can see from Figure 6 that $\gamma_4=\frac{1}{3}$ so we have
\begin{center}
$y_3=x^\frac{1}{3}(c_4+y_4)$
\end{center}
where $c_4$ is the next coefficient in the series and $y_4$ is the rest of the series. Substituting $y_3$ into $f_3$ gives
\begin{align*}
0=&\frac{2}{3}x^\frac{16}{3}(c_4+y_4)^3-x^3(c_4+y_4)^3+\frac{64}{3}x^3(c_4+y_4)-\frac{3073}{48}x^\frac{7}{3}(c_4+y_4) \\
&-32x^\frac{11}{3}(c_4+y_4)^2-6x^\frac{5}{3}(c_4+y_4)^2+\frac{16}{3}x^\frac{13}{3}(c_4+y_4)^2-\frac{4}{3}x^\frac{14}{3}(c_4+y_4)^2 \\
&-\frac{1}{9}x^\frac{13}{3}(c_4+y_4)+2x^\frac{17}{3}(c_4+y_4)^4-\frac{1}{4}x^\frac{8}{3}(c_4+y_4)^2-12x^{1}{3}(c_4+y_4) \\
&+48x^3(c_4+y_4)^2+4x^4(c_4+y_4)^2-\frac{9353}{216}x+16x^\frac{13}{3}(c_4+y_4)^3 \\
&+8x^\frac{8}{3}(c_4+y_4)-\frac{16}{3}x^\frac{10}{3}(c_4+y_4)+\frac{196}{3}x^\frac{5}{3}(c_4+y_4)+\frac{1}{216}x^\frac{14}{3}(c_4+y_4) \\
&-\frac{16}{3}x^5(c_4+y_4)^3+\frac{1}{12}x^5(c_4+y_4)^2-\frac{55}{27}x^\frac{11}{3}(c_4+y_4)+\frac{8}{9}x^4(c_4+y_4) \\
&-x^\frac{4}{3}(c_4+y_4)-\frac{1}{324}x^4-\frac{9217}{1728}x^2+\frac{14}{81}x^3+\frac{1}{27}x^\frac{11}{3}+\frac{16}{9}x^\frac{8}{3}-\frac{61}{324}x^\frac{10}{3}\\
&-\frac{119}{27}x^\frac{7}{3}+\frac{49}{9}x^\frac{4}{3}+\frac{104}{3}x^\frac{1}{3}+\frac{64}{3}x^\frac{5}{3}+\frac{1}{10368}x^\frac{13}{3} \tag{14} \\
\end{align*}
From the Newton Polygon (Figure 6), $\beta_4=\frac{1}{3}$, so we let the $x^\frac{1}{3}=0$ and solve for $c_4$.
\begin{center}
$\frac{104}{3}x^\frac{1}{3}-12c_4x^\frac{1}{3}=0$

so, $c_4=\frac{26}{9}$
\end{center}
Substituting $c_4$ into (14) and dividing by $x^\frac{1}{3}$
\begin{align*}
x^{-\frac{1}{3}}f_4&(x,x^\frac{1}{3}(c_4+y_4))=\frac{23}{9}x-\frac{4}{3}x^\frac{13}{3}y_4^2-\frac{5047}{27}x^\frac{10}{3}y_4+\frac{35152}{2187}x^5\\
&+\frac{40119049}{93312}x^4-\frac{40901}{216}x^2-\frac{5053}{324}x^3+144x^4y_4^2+\frac{11647}{324}x^\frac{11}{3}\\
&+\frac{319510}{729}x^\frac{8}{3}-\frac{66317}{243}x^\frac{10}{3}+\frac{1847}{81}x^\frac{7}{3}+160x^\frac{4}{3}-\frac{9353}{216}x^\frac{2}{3}+\frac{52}{9}x^5y_4^2\\
&+\frac{1352}{81}x^5y_4-\frac{16}{3}x^3y_4+\frac{11645}{27}x^4y_4+\frac{2}{3}x^5y_4^3-\frac{16}{3}x^\frac{14}{3}y_4^3+2x^\frac{16}{3}y_4^4\\
&-xy_4-\frac{3073}{48}x^2y_4+16x^4y_4^3-6x^\frac{4}{3}y_4^2-\frac{1}{4}x^\frac{7}{3}y_4^2-32x^\frac{10}{3}y_4^2\\
&-\frac{279695}{2187}x^\frac{14}{3}+\frac{2704}{27}x^\frac{16}{3}y_4^2+\frac{140608}{729}x^\frac{16}{3}y_4+24x^\frac{11}{3}y_4+\frac{92}{3}x^\frac{4}{3}y_4\\
&-\frac{1663}{216}x^\frac{13}{3}y_4-12y_4-\frac{9217}{1728}x^\frac{5}{3}-\frac{3601}{324}x^\frac{13}{3}+\frac{59}{9}x^\frac{7}{3}y_4+\frac{208}{9}x^\frac{16}{3}y_4^3 \\
&+\frac{913952}{6561}x^\frac{16}{3}+\frac{118}{3}x^\frac{8}{3}y_4^2+\frac{7388}{27}x^\frac{8}{3}y_4-\frac{10777}{81}x^\frac{14}{3}y_4-\frac{1661}{36}x^\frac{14}{3}y_4^2\\
&+4x^\frac{11}{3}y_4^2-x^\frac{8}{3}y_4^3 \tag{15} \\
\end{align*}
\begin{figure} [ht]
\centering
\caption{Newton Polygon for $f_4$}
\scalebox{.7}{\includegraphics{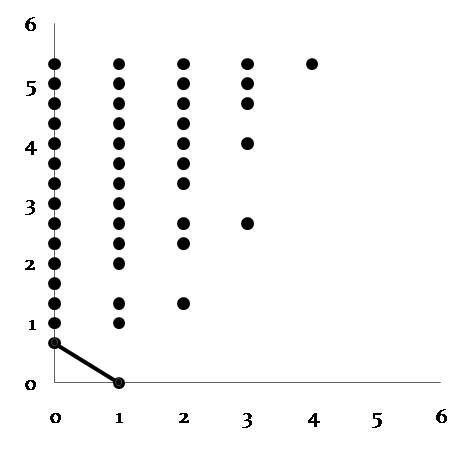}}
\end{figure}
\\
It can be seen from the new Newton Polygon (Figure 7) that $\gamma_5=\frac{2}{3}$, so 
\begin{center}
$y_4=x^\frac{2}{3}(c_5+y_5)$
\end{center}
\nopagebreak
Substituting $y_4$ into $f_4$ we find that
\begin{align*}
0=&\frac{23}{9}x+4x^5(c_5+y_5)^2+\frac{92}{3}x^2(c_5+y_5)+\frac{35152}{2187}x^5+\frac{40119049}{93312}x^4 \\
&-\frac{40901}{216}x^2-\frac{5053}{324}x^3-12x^\frac{2}{3}(c_5+y_5)+\frac{11647}{324}x^\frac{11}{3}+\frac{319510}{729}x^\frac{8}{3} \\
&-\frac{66317}{243}x^\frac{10}{3}+\frac{1847}{81}x^\frac{7}{3}+160x^\frac{4}{3}-\frac{9353}{216}x^\frac{2}{3}-\frac{279695}{2187}x^\frac{14}{3} \\
&+144x^\frac{16}{3}(c_5+y_5)^2-6x^\frac{8}{3}(c_5+y_5)^2-x^\frac{14}{3}(c_5+y_5)^3-\frac{1663}{216}x^5(c_5+y_5) \\
&-\frac{16}{3}x^\frac{11}{3}(c_5+y_5)-\frac{1}{4}x^\frac{11}{3}(c_5+y_5)^2+24x^\frac{13}{3}(c_5+y_5) \\
&-\frac{10777}{81}x^\frac{16}{3}(c_5+y_5)+\frac{2}{3}x^7(c_5+y_5)^3+\frac{1352}{81}x^\frac{17}{3}(c_5+y_5) \\
&+\frac{140608}{729}x^6(c_5+y_5)-32x^\frac{14}{3}(c_5+y_5)^2-\frac{3073}{48}x^\frac{8}{3}(c_5+y_5) \\
&+\frac{59}{9}x^3(c_5+y_5)+\frac{52}{9}x^\frac{19}{3}(c_5+y_5)^2+\frac{208}{9}x^\frac{22}{3}(c_5+y_5)^3-x^\frac{5}{3}(c_5+y_5) \\
&+2x^8(c_5+y_5)^4+\frac{118}{3}x^4(c_5+y_5)^2-\frac{1661}{36}x^6(c_5+y_5)^2+16x^6(c_5+y_5)^3 \\
&+\frac{2704}{27}x^\frac{20}{3}(c_5+y_5)^2+\frac{11645}{27}x^\frac{14}{3}(c_5+y_5)-\frac{5047}{27}x^4(c_5+y_5) \\
&+\frac{7388}{27}x^\frac{10}{3}(c_5+y_5)-\frac{16}{3}x^\frac{20}{3}(c_5+y_5)^3-\frac{4}{3}x^\frac{17}{3}(c_5+y_5)^2-\frac{9217}{1728}x^\frac{5}{3} \\
&-\frac{3601}{324}x^\frac{13}{3}+\frac{913952}{6561}x^\frac{16}{3} \tag{16} \\
\end{align*}
We can see from the Newton Polygon (Figure 7) that $\beta_5=\frac{2}{3}$, so we can set the $x^\frac{2}{3}$ terms equal zero and solve for $c_5$.
\begin{center}
$\frac{9353}{216}x^\frac{2}{3}-12c_5x^\frac{2}{3}=0$

so, $c_5=\frac{9353}{2592}$
\end{center}
This process could be continued to calculate more terms, but we will stop here. So, an explicit solution for $f(x,y)$ is
\begin{equation*}
y=\frac{9353}{2592}x^\frac{10}{3}+\frac{26}{9}x^\frac{8}{3}+\frac{1}{12}x^\frac{7}{3}+\frac{2}{3}x^2+2x^\frac{4}{3}
\end{equation*}
 
We must now examine the case where the solution to $c_1$ is imaginary. There are two imaginary solutions, but we will only show the process for one of these cases because they are quite similar. Let $c_1=-1+\sqrt{2}i$. We substitute $c_1$ into (5) and divide by $x^4$.
\begin{align*}
x^{-4}f_1(x,x^2(c_1+y_1))&=f_1(x,y_1) \\
&=x+4i\sqrt{3}x^\frac{2}{3}-3i\sqrt{3}y_1^2-16x^\frac{4}{3}-24x^\frac{4}{3}y_1^2-8x^\frac{4}{3}y_1^3 \\
&+2x^\frac{4}{3}y_1^4+64x^\frac{4}{3}y_1-24ix^\frac{4}{3}3^\frac{1}{2}y_1^2+6i\sqrt{3}y_1-4i\sqrt{3}x^\frac{2}{3}y_1 \\
&+4x^\frac{2}{3}y_1-2x^\frac{2}{3}y_1^2+16i\sqrt{3}x^\frac{4}{3}+8i\sqrt{3}x^\frac{4}{3}y_1^3+6y_1+3y_1^2 \\
&-y_1^3+4x^\frac{2}{3} \tag{17}\\
\end{align*}
\begin{figure} [ht]
\centering
\caption{Newton Polygon for $f_1$}
\scalebox{.7}{\includegraphics{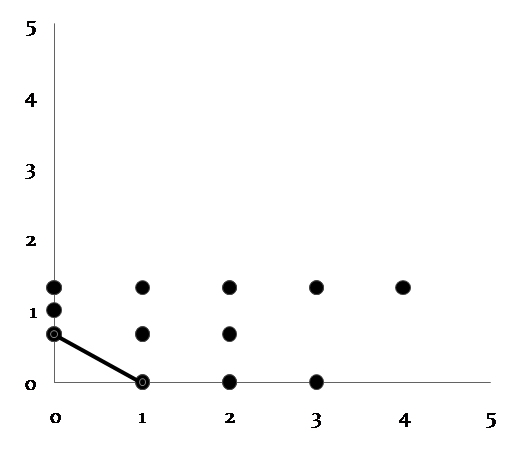}}
\end{figure}
\\
We can see from the new Newton Polygon (Figure 8) that $\gamma_2=\frac{2}{3}$, so 
\begin{center}
$y_1=x^\frac{2}{3}(c_2+y_2)$
\end{center}
where $c_2$ is the second coefficient in the series and 
\begin{center}
$y_2=c_3x^{\gamma_3}+c_4x^{\gamma_3+\gamma_4}+c_5x^{\gamma_3+\gamma_4+\gamma_5}+\cdots$
\end{center}
We substitute $y_1$ into $f_1$ to get
\begin{align*}
0=&x+4i\sqrt{3}x^\frac{2}{3}-4i\sqrt{3}x^\frac{4}{3}(c_2+y_2)-16x^\frac{4}{3}-24x^\frac{8}{3}(c_2+y_2)^2-8x^\frac{10}{3}(c_2+y_2)^3\\
 &+2x^4(c_2+y_2)^4+64x^2(c_2+y_2)+6i\sqrt{3}x^\frac{2}{3}(c_2+y_2)-3i\sqrt{3}x^\frac{4}{3}(c_2+y_2)^2 \\
&-24i\sqrt{3}x^\frac{8}{3}(c_2+y_2)^2+4x^\frac{4}{3}(c_2+y_2)-2x^2(c_2+y_2)^2+16i\sqrt{3}x^\frac{4}{3} \\
&+8i\sqrt{3}x^\frac{10}{3}(c_2+y_2)^3+6x^\frac{2}{3}(c_2+y_2)+3x^\frac{4}{3}(c_2+y_2)^2-x^2(c_2+y_2)^3+4x^\frac{2}{3} \tag{18} \\ 
\end{align*}
\pagebreak
\\
From the Newton Polygon (Figure 8), $\beta_2=\frac{2}{3}$, so we can let the $x^\frac{2}{3}$ terms equal zero and solve for $c_2$.
\begin{center}
$6c_2x^\frac{2}{3}+6i\sqrt{3}c_2x^\frac{2}{3}+4i\sqrt{3}x^\frac{2}{3}+4x^\frac{2}{3}=0$

so, $c_2=-\frac{2}{3}$
\end{center}
Now, substituting $c_2$ into (18) and dividing by $x^\frac{2}{3}$ gives
\begin{align*}
x^{-\frac{2}{3}}f_2&(x,x^\frac{2}{3}(c_2+y_2))=-\frac{32}{3}x^2+\frac{64}{27}x^\frac{8}{3}+\frac{32}{81}x^\frac{10}{3}-\frac{1168}{27}x^\frac{4}{3}+x^\frac{1}{3}-\frac{52}{3}x^\frac{2}{3}\\
&+6y_2+32x^2y_2-24x^2y_2^2-\frac{16}{3}x^\frac{10}{3}y_2^3-\frac{64}{27}x^\frac{10}{3}y_2+\frac{16}{3}x^\frac{10}{3}y_2^2 \\
&-\frac{32}{3}x^\frac{8}{3}y_2+16x^\frac{8}{3}y_2^2+\frac{196}{3}x^\frac{4}{3}y_2+3x^\frac{2}{3}y_2^2-x^\frac{4}{3}y_2^3-8x^\frac{8}{3}y_2^3\\
&+2x^\frac{10}{3}y_2^4-3i\sqrt{3}x^\frac{2}{3}y_2^2-24i\sqrt{3}x^2y_2^2+32i\sqrt{3}x^2y_2+\frac{32}{3}i\sqrt{3}x^\frac{8}{3}y_2\\
&+8i\sqrt{3}x^\frac{8}{3}y_2^3-16i\sqrt{3}x^\frac{8}{3}y_2^2+\frac{52}{3}i\sqrt{3}x^\frac{2}{3}-\frac{64}{27}i\sqrt{3}x^\frac{8}{3}-\frac{32}{3}i\sqrt{3}x^2\\
&+6i\sqrt{3}y_2 \tag{19} \\
\end{align*}
\begin{figure} [ht]
\centering
\caption{Newton Polygon for $f_2$}
\scalebox{.7}{\includegraphics{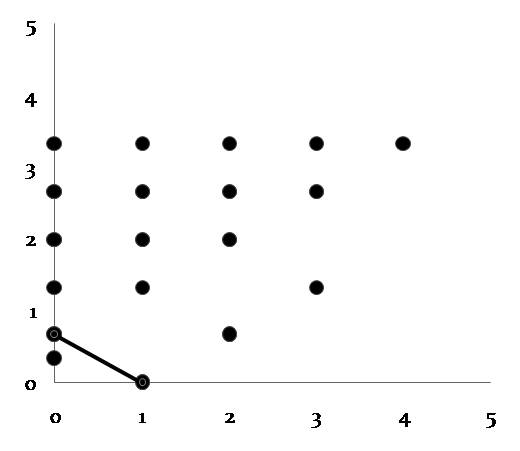}}
\end{figure}
\\
From Figure 9 we see that $\gamma_3=\frac{1}{3}$, so
\begin{center}
$y_2=x^\frac{1}{3}(c_3+y_3)$
\end{center}
where $c_3$ is the next coefficient in the series and $y_3$ is the rest of the series. Substituting $y_2$ into $f_2$ gives
\begin{align*}
0=&\frac{16}{3}x^4(c_3+y_3)^2-\frac{32}{3}x^2+\frac{64}{27}x^\frac{8}{3}+\frac{32}{81}x^\frac{10}{3}-\frac{1168}{27}x^\frac{4}{3}+x^\frac{1}{3} \\
&-\frac{52}{3}x^\frac{2}{3}+6x^\frac{1}{3}(c_3+y_3)-x^\frac{7}{3}(c_3+y_3)^3-\frac{64}{27}x^\frac{11}{3}(c_3+y_3)-\frac{32}{3}x^3(c_3+y_3) \\
&-24x^\frac{8}{3}(c_3+y_3)^2+2x^\frac{14}{3}(c_3+y_3)^4-8x^\frac{11}{3}(c_3+y_3)^3+\frac{196}{3}x^\frac{5}{3}(c_3+y_3) \\
&+16x^\frac{10}{3}(c_3+y_3)^2+32x^\frac{7}{3}(c_3+y_3)+3x^\frac{4}{3}(c_3+y_3)^2-\frac{16}{3}x^\frac{13}{3}(c_3+y_3)^3 \\
&+6i\sqrt{3}x^\frac{1}{3}(c_3+y_3)+\frac{52}{3}i\sqrt{3}x^\frac{2}{3}-\frac{64}{27}i\sqrt{3}x^\frac{8}{3}-\frac{32}{3}i\sqrt{3}x^2 \\
&+\frac{32}{3}i\sqrt{3}x^3(c_3+y_3)+8i\sqrt{3}x^\frac{11}{3}(c_3+y_3)^3-16i\sqrt{3}x^\frac{10}{3}(c_3+y_3)^2 \\
&-24i\sqrt{3}x^\frac{8}{3}(c_3+y_3)^2-3i\sqrt{3}x^\frac{4}{3}(c_3+y_3)^2+32i\sqrt{3}x^\frac{7}{3}(c_3+y_3) \tag{20} \\
\end{align*}
We can see from the Newton Polygon (Figure 9) that $\beta_3=\frac{1}{3}$, so we can let the $x^\frac{1}{3}$ terms equal zero and solve for $c_3$.
\begin{center}
$6c_3x^\frac{1}{3}+1x^\frac{1}{3}+6i\sqrt{3}c_3x^\frac{1}{3}=0$

so, $c_3=\frac{1}{6(\sqrt{3}+1)}$
\end{center}
\nopagebreak
Substituting $c_3$ into (20) and dividing by $x^\frac{1}{3}$ gives us
\begin{align*}
x^{-\frac{1}{3}}f_3&(x,x^\frac{1}{3}(c_3+y_3))=-\frac{128i\sqrt{3}x^\frac{11}{3}y_3^2}{3(i\sqrt{3}+1)^4}+\frac{8i\sqrt{3}x^2y_3^3}{(i\sqrt{3}+1)^4}-\frac{64i\sqrt{3}x^\frac{7}{3}y_3}{(i\sqrt{3}+1)^4} \\
&-\frac{8i\sqrt{3}y_3x}{(i\sqrt{3}+1)^4}-\frac{128x^\frac{8}{3}}{9(i\sqrt{3}+1)^4}-\frac{244x^\frac{10}{3}}{81(i\sqrt{3}+1)^4}-\frac{1904x^\frac{7}{3}}{27(i\sqrt{3}+1)^4} \\
&-\frac{1568x^\frac{4}{3}y_3}{3(i\sqrt{3}+1)^4}-\frac{128x^\frac{11}{3}y_3^2}{3(i\sqrt{3}+1)^4}+\frac{128x^\frac{11}{3}y_3}{9(i\sqrt{3}+1)^4}-\frac{96xy_3^2}{(i\sqrt{3}+1)^4} \\
&+\frac{3073x^2y_3}{6(i\sqrt{3}+1)^4}+\frac{96y_3}{(i\sqrt{3}+1)^4}+\frac{8xy_4}{(i\sqrt{3}+1)^4}-\frac{4x^2y_3^2}{(i\sqrt{3}+1)^4} \\
&+\frac{8x^2y_3^3}{(i\sqrt{3}+1)^4}-\frac{64x^\frac{7}{3}y_3}{(i\sqrt{3}+1)^4}-\frac{8x^\frac{11}{3}}{27(i\sqrt{3}+1)^4}-\frac{8i\sqrt{3}x^4y_3}{9(i\sqrt{3}+1)^4} \\
&+\frac{128i\sqrt{3}x^4y_3^3}{3(i\sqrt{3}+1)^4}+\frac{x^\frac{13}{3}}{648(i\sqrt{3}+1)^4}+\frac{2x^4}{81(i\sqrt{3}+1)^4}+\frac{784x^\frac{4}{3}}{9(i\sqrt{3}+1)^4} \\
&-\frac{512x^3y_3^2}{(i\sqrt{3}+1)^4}+\frac{440x^\frac{10}{3}y_3}{27(i\sqrt{3}+1)^4}-\frac{384x^\frac{7}{3}y_3^2}{(i\sqrt{3}+1)^4}+\frac{128x^3y_3}{3(i\sqrt{3}+1)^4} \\
&-\frac{64x^4y_3^2}{3(i\sqrt{3}+1)^4}+\frac{128x^4y_3^3}{3(i\sqrt{3}+1)^4}-\frac{32x^\frac{10}{3}y_3^2}{(i\sqrt{3}+1)^4}+\frac{256x^\frac{10}{3}y_3^3}{(i\sqrt{3}+1)^4} \\
&+\frac{8x^4y_3}{9(i\sqrt{3}+1)^4}-\frac{x^\frac{13}{3}}{27(i\sqrt{3}+1)^4}-\frac{16x^\frac{13}{3}}{i\sqrt{3}+1)^4}+\frac{32x^\frac{13}{3}y_3^3}{3(i\sqrt{3}+1)^4} \\
&-\frac{2x^\frac{13}{3}y_3^2}{3(i\sqrt{3}+1)^4}+\frac{2i\sqrt{3}x^4}{81(i\sqrt{3}+1)^4}+\frac{128i\sqrt{3}x^\frac{8}{3}}{9(i\sqrt{3}+1)^4}+\frac{9217i\sqrt{3}x^2}{216(i\sqrt{3}+1)^4} \\
&+\frac{512i\sqrt{3}x^\frac{5}{3}}{3(i\sqrt{3}+1)^4}+\frac{9353ix\sqrt{3}}{27(i\sqrt{3}+1)^4}-\frac{112i\sqrt{3}x^3}{81(i\sqrt{3}+1)^4}+\frac{8i\sqrt{3}x^\frac{11}{3}}{27(i\sqrt{3}+1)^4} \\
&-\frac{96i\sqrt{3}y_3}{(i\sqrt{3}+1)^4}+\frac{2i\sqrt{3}x^\frac{13}{3}y_3^2}{(i\sqrt{3}+1)^4}-\frac{i\sqrt{3}x^\frac{13}{3}y_3}{27(i\sqrt{3}+1)^4}-\frac{16i\sqrt{3}x^\frac{13}{3}y_3^4}{(i\sqrt{3}+1)^4} \\
&-\frac{3073i\sqrt{3}x^2y_3}{6(i\sqrt{3}+1)^4}-\frac{1568i\sqrt{3}x^\frac{4}{3}}{3(i\sqrt{3}+1)^4}+\frac{440i\sqrt{3}x^\frac{10}{3}y_3}{27(i\sqrt{3}+1)^4}+\frac{32i\sqrt{3}x^\frac{10}{3}y_3}{(i\sqrt{3}+1)^4} \\
&-\frac{128i\sqrt{3}x^3y_3}{3(i\sqrt{3}+1)^4}-\frac{112x^3}{81(i\sqrt{3}+1)^4}+\frac{1024x^\frac{8}{3}y_3}{3(i\sqrt{3}+1)^4}+\frac{9217x^2}{216(i\sqrt{3}+1)^4} \\
&+\frac{1664x^\frac{1}{3}}{3(i\sqrt{3}+1)^4}+\frac{9353x}{27(i\sqrt{3}+1)^4}-\frac{512x^\frac{5}{3}}{3(i\sqrt{3}+1)^4}+\frac{384i\sqrt{3}x^\frac{7}{3}y_3}{(i\sqrt{3}+1)^4} \tag{21} \\
\end{align*}
\pagebreak
\\
By continuing this process it can be found that $c_4=\frac{52}{9(-1+i\sqrt{3})}$ and\\ $c_5=\frac{9353}{1296(1+i\sqrt{3})}$. So, a second explicit solution for $f(x,y)$ is
\begin{equation*}
y=-1+i\sqrt{3}x^\frac{4}{3}-\frac{2}{3}x^2-\frac{x^\frac{7}{3}}{6(i\sqrt{3}+1)}+\frac{52x^\frac{8}{3}}{9(-1+i\sqrt{3})}+\frac{9353x^\frac{10}{3}}{1296(i\sqrt{3}+1)}
\end{equation*}
Using a similar process, the explicit solution of $f(x,y)$ in the case where\\ $c_1=-1-i\sqrt{3}$ can be found. So, we have that the explicit solutions for $f(x,y)$ are
\begin{align*}
&y=\frac{9353}{2592}x^\frac{10}{3}+\frac{26}{9}x^\frac{8}{3}+\frac{1}{12}x^\frac{7}{3}+\frac{2}{3}x^2+2x^\frac{4}{3}+\cdots \\
&y=-1+i\sqrt{3}x^\frac{4}{3}-\frac{2}{3}x^2-\frac{x^\frac{7}{3}}{(6(i\sqrt{3}+1)}+\frac{52x^\frac{8}{3}}{9(-1+i\sqrt{3})}+\frac{9353x^\frac{10}{3}}{1296(i\sqrt{3}+1)}+\cdots \\
&y=(-1-i\sqrt{3})x^\frac{4}{3}-\frac{2}{3}x^2+\frac{x^\frac{7}{3}}{(6(-1+i\sqrt{3})}-\frac{52x^\frac{8}{3}}{9(i\sqrt{3}+1)}-\frac{9353x^{10}{3}}{1296(-1+i\sqrt{3})}+\cdots \\
\end{align*}
\bigskip
\emph{Example 3:} \\
\begin{figure} [ht]
\centering
\caption{Newton Polygon for $f(x,y)$}
\scalebox{.7}{\includegraphics{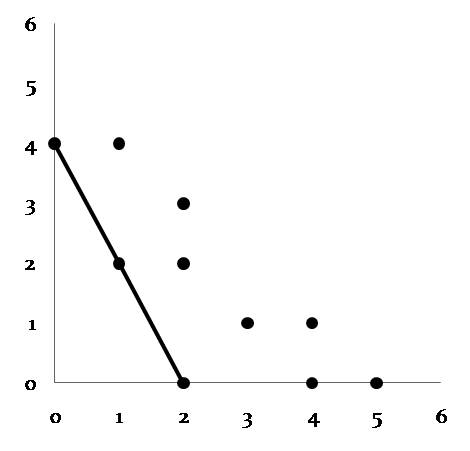}}
\end{figure}

In this example, we find that multiple terms of a Puiseux expansion can be calculated before the Puiseux jets differ.  Here, the jets will appear the same for the first three terms and then split at $\frac{11}{2}$.
\medskip
\begin{align*}
f(x,y)=y^2+2x^2y+x^4+x^2y^2+xy^3+\frac{1}{4}y^4+x^4y+x^3y^2-\frac{1}{2}xy^4-\frac{1}{2}y^5 \tag{1} \\
\end{align*}
\medskip

We know that this polynomial has solutions of the form 
\begin{align*}
y=c_1x^{\gamma_1}+c_2x^{\gamma_1+\gamma_2}+c_3x^{\gamma_1+\gamma_2+\gamma_3}+\cdots. \tag{2} \\
\end{align*}
The possible choices for $\gamma_1$ correspond to the negatives of the slopes of the lower segments of the Newton Polygon (Figure 10). The slopes of the lower segments of the Newton Polygon corresponding to (1) are both -2. We will start by considering this case where $\gamma_1=-(-2)=2$.

Let $\gamma_1=2$ Then by factoring $x^{\gamma_1}$ (or $x^2$) out of (2) we obtain, 
\begin{align*}
y=x^2(c_1+y_1) \tag{3} \\
\end{align*}
where $c_1$ is the first coefficient and $y_1$ is
\begin{align*}
y_1=c_2x^{\gamma_2}+c_3x^{\gamma_2+\gamma_3}+c_4x^{\gamma_2+\gamma_3+\gamma_4}+\cdots \tag{4} \\
\end{align*}
\\ i.e. the rest of the series.
\\ So, substituting (3) into $f$ we get, 
\begin{align*}
0=&x^4(c_1+y_1)^2+2x^4(c_1+y_1)+x^4+x^6(c_1+y_1)^2+x^7(c_1+y_1)^3\\
&+\frac{1}{4}x^8(c_1+y_1)^4+x^6(c_1+y_1)+x^7(c_1+y_1)^2-\frac{1}{2}x^9(c_1+y_1)^4\\
&-\frac{1}{2}x^{10}(c_1+y_1)^5 \tag{5} \\
\end{align*}

The vertical intercept on the Newton Polygon (Figure 10), which we call $\beta$, gives the terms of lowest degree in $x$ alone. These terms must sum to zero because $f(x,y)=0$. In this case, $\beta=4$ so, $x^4c_1^2+x^4c_1+x^4=0$.  The $x^n$ part of the terms can be divided out and will not be included in subsequent steps. Solving for $c_1$ we obtain $c_1=-1$ 
Substituting  $c_1$ into (5) and dividing by $x^4$ gives, 

\begin{align*}
x^{-4}&f(x,x^2(c_1+y_1))=f_1(x,y_1)=\\
&=-\frac{1}{2}x^5+\frac{1}{2}x^6+\frac{1}{4}x^4+y_1^2-x^2y_1+x^2y_1^2+x^3y_1\\
&-2x^3y_1^2+x^3y_1^3-x^4y_1+\frac{3}{2}x^4y_1^2-x^4y_1^3+\frac{1}{4}x^4y_1^4+2x^5y_1-3x^5y_1^2\\
&+2x^5y_1^3-\frac{1}{2}x^5y_1^4-\frac{5}{2}x^6y_1+5x^6y_1^2-5x^6y_1^3+\frac{5}{2}x^6y_1^4-\frac{1}{2}x^6y_1^5 \tag{6} \\
 \end{align*}
\begin{figure} [ht]
\centering
\caption{Newton Polygon for $f_1$}
\scalebox{.7}{\includegraphics{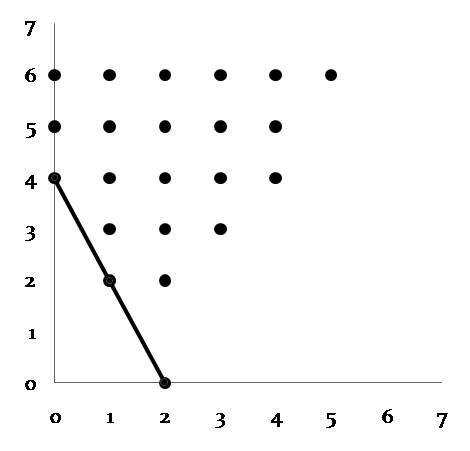}}
\end{figure}
\\
We call this new polynomial $f_1(x,y_1)$, so now that we have our new function $f_1(x,y_1)$, and our new goal is to find $c_2$. We find the new Newton Polygon (Figure 11), from which we get $\beta_2=4$ and $\gamma_2=2$.  Similar to (3), we then have

\begin{align*}
 y_1=x^2(c_2+y_2) \tag{7} \\
\end{align*}

 where $c_2$ is the coefficient of the second term in the series and $y_2$ is

\begin{align*}
 y_2=c_3x^{\gamma_3}+c_4x^{\gamma_3+\gamma_4}+\cdots \tag{8} \\
\end{align*}  

We substitute (7) into $f_1$ to get
\begin{align*}
0=&-\frac{1}{2}x^5+\frac{1}{2}x^6+\frac{1}{4}x^4+x^4(c_2+y_2)^2-x^4(c_2+y_2)+x^6(c_2+y_2)^2\\
&+x^5(c_2+y_2)-2x^7(_c2+y_2)^2+x^9(c_2+y_2)^3-x^6(c_2+y_2)\\
&+\frac{3}{2}x^8(c_2+y_2)^2-x^{10}(c_2+y_2)^3+\frac{1}{4}x^{12}(c_2+y_2)^4+2x^7(c_2+y_2)\\
&-3x^9(c_2+y_2)^2+2x^{11}(c_2+y_2)^3-\frac{1}{2}x^{13}(c_2+y_2)^4-\frac{5}{2}x^8(c_2+y_2)\\
&+5x^{10}(c_2+y_2)^2-5x^{12}(c_2+y_2)^3+\frac{5}{2}x^{14}(c_2+y_2)^4-\frac{1}{2}x^{16}(c_2+y_2)^5 \tag{9}  \\
\end{align*}
\pagebreak
Again, lowest terms must cancel so, since $\beta_2=4$, we see that $\frac{1}{4}+c_2^2-c_2=0$.  Solving for $c_2$ gives us $c_2=\frac{1}{2}$, and we substitute this $c_2$ into (9) and divide by $x^{-\beta}$ giving

\begin{align*}
x^{-4}&f(x,x^2(c_2+y_2))=f_2(x,y_2)=\\
&=-\frac{5}{8}x^5+\frac{9}{8}x^6+\frac{1}{4}x^7-\frac{1}{32}x^9-\frac{39}{64}x^8-\frac{1}{64}x^{12}\\
&+\frac{5}{32}x^{10}+\frac{1}{4}x^2-\frac{7}{8}x^4+\frac{1}{2}x^3+y_2^2+x^2y_2^2+xy_2-2x^3y_2^2-\frac{9}{4}x^5y_2\\
&-\frac{3}{2}x^5y_2^2+x^5y_2^3-x^4y_2+\frac{3}{2}x^4y_2^2+\frac{17}{4}x^6y_2+\frac{7}{2}x^6y_2^2-x^6y_2^3-\frac{29}{8}x^8y_2\\
&-\frac{57}{8}x^8y_2^2-\frac{9}{2}x^8y_2^3+\frac{1}{4}x^8y_2^4+\frac{3}{2}x^7y_2+3x^7y_2^2+2x^7y_2^3-\frac{1}{4}x^9y_2\\
&-\frac{3}{4}x^9y_2^2-x^9y_2^3-\frac{1}{2}x^9y_2^4+\frac{5}{4}x^{10}y_2+\frac{15}{4}x^{10}y_2^2+5x^{10}y_2^3+\frac{5}{2}x^{10}y_2^4\\
&-\frac{5}{32}x^{12}y_2-\frac{5}{8}x^{12}y_2^2-\frac{5}{4}x^{12}y_2^3-\frac{5}{4}x^{12}y_2^4-\frac{1}{2}x^{12}y_2^5 
\tag{10} \\
\end{align*}
\begin{figure} [ht]
\centering
\caption{Newton Polygon for $f_2$}
\scalebox{.7}{\includegraphics{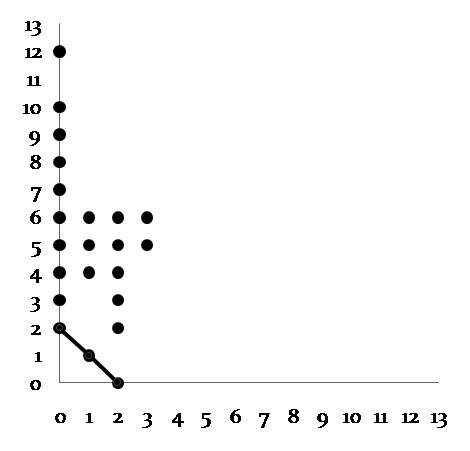}}
\end{figure}
\\
We must now find $c_3$.  From the Newton Polygon for $f_2$ (Figure 12), we can see that $\beta_3=2$ and $\gamma_3=1$.  Therefore,

\begin{align*}
y_2=x(c_3+y_3) \tag{11} \\
\end{align*}
where $c_3$ is the coefficient of the third term in the series and $y_3$ is
\begin{align*}
y_3=c_4x^{\gamma_4}+c_5x^{\gamma_4+\gamma_5}+\cdots \tag{12} \\
\end{align*}

Substituting (11) into $f_2$ we get

\begin{align*}
0=&\frac{15}{4}x^{12}(c_3+y_3)^2-\frac{57}{8}x^{10}(c_3+y_3)^2-{29}{8}x^9(c_3+_y3)-x^{12}(c_3+y_3)^3\\
&-2x^5(c_3+y_3)^2-\frac{3}{4}x^{11}(c_3+y_3)^2-x^9(c_3+y_3)^3+x^2(c_3+y_3)\\
&+\frac{3}{2}x^8(c_3+y_3)-\frac{5}{4}x^{16}(c_3+y_3)^4-\frac{3}{2}x^7(c_3+y_3)^2+\frac{7}{2}x^8(c_3+y_3)^2\\
&-\frac{1}{4}x^{10}(c_3+y_3)+\frac{5}{2}x^{14}(c_3+y_3)^4-\frac{5}{32}x^{13}(c_3+y_3)-\frac{1}{2}x^{13}(c_3+y_3)^4\\
&+\frac{5}{4}x^{11}(c_3+y_3)+3x^9(c_3+y_3)^2+\frac{1}{4}x^{12}(c_3+y_3)^4+x^8(c_3+y_3)^3\\
&-x^5(c_3+y_3)+\frac{3}{2}x^6(c_3+y_3)^2-\frac{1}{2}x^{17}(c_3+y_3)^5-\frac{5}{8}x^{14}(c_3+y_3)^2\\
&-\frac{9}{2}x^{11}(c_3+y_3)^3+\frac{1}{4}x^2-\frac{7}{8}x^4+\frac{1}{2}x^3-\frac{9}{4}x^6(c_3+y_3)-\frac{5}{4}x^{15}(c_3+y_3)^3\\
&+2x^{10}(c_3+y_3)^3-\frac{5}{8}x^5+\frac{9}{8}x^6+\frac{1}{4}x^7-\frac{1}{32}x^9-\frac{39}{64}x^8+\frac{17}{4}x^7(c_3+y_3)\\
&-\frac{1}{64}x^{12}+\frac{5}{32}x^{10}+x^4(c_3+y_3)^2+x^2(c_3+y_3)^2+5x^{13}(c_3+y_3)^3 \tag{13} \\
\end{align*}

Since $\beta_3=2$ (which can be seen from Figure 12), when lowest terms cancel, we find that $c_3+\frac{1}{4}+c_3^2=0$ so $c_3=-\frac{1}{2}$.  We substitute this value for $c_3$ into (12) and divide by $x^{-\beta}$ to get
\begin{align*}
x^{-2}&f_3(x,x(c_3+y_3))=f_3(x,y_3)=\\
&=\frac{1}{2}x-{5}{8}x^2+\frac{21}{8}x^4-\frac{5}{8}x^3-\frac{9}{4}x^5-\frac{39}{64}x^6+\frac{85}{32}x^7-\frac{1}{4}x^9\\
&-\frac{7}{4}x^8+\frac{17}{16}x^{10}-\frac{37}{64}x^{11}+\frac{1}{64}x^{15}-\frac{5}{64}x^{14}+\frac{5}{32}x^{13}+y_3^2-\frac{5}{8}x^{12}y_3\\
&+\frac{25}{8}x^{12}y_3^2-\frac{37}{8}x^{10}y_3+\frac{45}{8}x^{10}y_3^2-\frac{11}{8}x^9y_3+\frac{23}{4}x^5y_3-\frac{3}{2}x^5y_3^2+6x^9y_3^2\\
&-\frac{9}{2}x^9y_3^3+\frac{67}{8}x^8y_3-\frac{59}{8}x^7y_3+\frac{9}{2}x^7y_3^2-\frac{81}{8}x^8y_3^2+\frac{5}{8}x^{14}y_3-\frac{15}{8}x^{14}y_3^2\\
&+\frac{5}{2}x^{14}y_3^3-\frac{5}{4}x^{14}y_3^4+\frac{123}{32}x^{11}y_3-5x^{12}y_3^3+\frac{5}{2}x^{12}y_3^4+2x^8y_3^3-\frac{5}{4}x^6y_3\\
&+2x^6y_3^2-\frac{33}{4}x^{11}y_3^2+6x^{11}y_3^3-\frac{3}{2}x^{10}y_3^3-\frac{15}{4}x^4y_3+\frac{3}{2}x^4y_3^2-x^2y_3+x^2y_3^2\\
&-\frac{15}{16}x^{13}y_3+\frac{15}{8}x^{13}y_3^2-\frac{5}{4}x^{13}y_3^3+x^3y_3-2x^3y_3^2-x^7y_3^3-\frac{1}{2}x^{11}y_3^4\\
&+\frac{1}{4}x^{10}y_3^4+x^6y_3^3-\frac{5}{32}x^{15}y_3+\frac{5}{8}x^{15}y_3^2-\frac{5}{4}x^{15}y_3^3+\frac{5}{4}x^{15}y_3^4-\frac{1}{2}x^{15}y_3^5 \tag{14} \\
\end{align*}
\begin{figure} [ht]
\centering
\caption{Newton Polygon for $f_3$}
\scalebox{.7}{\includegraphics{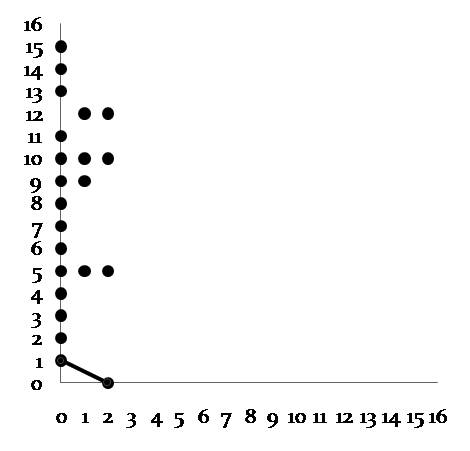}}
\end{figure}
\\
To find $c_4$, we repeat the process for $f_3$.  Looking at the Newton Polygon (Figure 13), we see that $\beta_4=1$ and $\gamma_4={\frac{1}{2}}$.  From this we have

\begin{align*}
y_3=x^{\frac{1}{2}}(c_4+y_4) \tag{15} \\
\end{align*}
where $c_4$ is the coefficient of the fourth term in the Puiseux series and $y_4$ is

\begin{align*}
y_4=c_4x^{\gamma_4}+c_5x^{\gamma_4+\gamma_5}+\cdots \tag{16} \\
\end{align*}
\pagebreak 
\\ Substituting (15) into $f_3$, we get

\begin{align*}
0=&\frac{3}{2}x^5(c_4+y_4)^2+\frac{123}{32}x^{\frac{23}{2}}(c_4+y_4)+\frac{1}{2}x+6x^{10}(c_4+y_4)^2-\frac{11}{8}x^{\frac{19}{2}}(c_4+y_4)\\
&-\frac{15}{16}x^{\frac{27}{2}}(c_4+y_4)-2x^4(c_4+y_4)^2+2x^{\frac{19}{2}}(c_4+y_4)^3-\frac{5}{32}x^{\frac{31}{2}}(c_4+y_4)\\&-\frac{5}{4}x^{\frac{33}{2}}(c_4+y_4)^3+x^{\frac{7}{2}}(c_4+y_4)-x^{\frac{17}{2}}(c_4+y_4)^3+\frac{5}{8}x^{16}(c_4+y_4)^2\\
&+\frac{5}{8}x^{\frac{29}{2}}(c_4+y_4)+\frac{9}{2}x^8(c_4+y_4)^2-\frac{1}{2}x^{13}(c_4+y_4)^4-\frac{1}{2}x^{\frac{35}{2}}(c_4+y_4)^5\\
&+\frac{45}{8}x^{11}(c_4+y_4)^2-\frac{5}{4}x^{\frac{13}{2}}(c_4+y_4)+\frac{1}{4}x^{12}(c_4+y_4)^4-5x^{\frac{27}{2}}(c_4+y_4)^3\\
&-\frac{59}{8}x^{\frac{15}{2}}(c_4+y_4)-\frac{15}{4}x^{\frac{9}{2}}(c_4+y_4)-\frac{81}{8}x^9(c_4+y_4)^2-x^{\frac{5}{2}}(c_4+y_4)\\
&-\frac{3}{2}x^{\frac{23}{2}}(c_4+y_4)^3+\frac{5}{2}x^{14}(c_4+y_4)^4-\frac{33}{4}x^{12}(c_4+y_4)^2-\frac{5}{4}x^{\frac{29}{2}}(c_4+y_4)^3\\
&+\frac{67}{8}x^{\frac{17}{2}}(c_4+y_4)+\frac{25}{8}x^{13}(c_4+y_4)^2+x^{\frac{15}{2}}(c_4+y_4)^3-\frac{3}{2}x^6(c_4+y_4)^2\\
&+\frac{15}{8}x^{14}(c_4+y_4)^2-\frac{5}{8}x^2+\frac{21}{8}x^4-\frac{5}{8}x^3+\frac{23}{4}x^{\frac{11}{2}}(c_4+y_4)\\
&-\frac{9}{2}x^{\frac{21}{2}}(c_4+y_4)^3-\frac{5}{4}x^{16}(c_4+y_4)^4-\frac{9}{4}x^5-\frac{39}{64}x^6+\frac{85}{32}x^7+\frac{5}{4}x^{17}(c_4+y_4)^4\\
&-\frac{37}{8}x^{\frac{21}{2}}(c_4+y_4)-\frac{1}{4}x^9-\frac{7}{4}x^8+\frac{5}{2}x^{\frac{31}{2}}(c_4+y_4)^3+\frac{17}{16}x^{10}\\
&-\frac{37}{64}x^{11}+\frac{1}{64}x^{15}-\frac{5}{64}x^{14}+\frac{5}{32}x^{13}+2x^7(c_4+y_4)^2-\frac{5}{8}x^{\frac{25}{2}}(c_4+y_4)\\
&+x^3(c_4+y_4)^2+6x^{\frac{25}{2}}(c_4+y_4)^3-\frac{15}{8}x^{15}(c_4+y_4)^2+x(c_4+y_4)^2 \tag{17} \\
\end{align*}

Since $\beta_4=1$, when lowest terms cancel, we find that $\frac{1}{2}+c_4^2=0$ so $c_4=\pm\frac{\sqrt{2}i}{2}$.  So the two solutions given by the Newton Puiseux Algorithm are
\begin{align*}
y=-x^2+\frac{x^4}{2}-\frac{x^5}{2}+\frac{(-2x)^{\frac{11}{2}}}{64}+\cdots \\
\end{align*}
\begin{align*}
y=-x^2+\frac{x^4}{2}-\frac{x^5}{2}-\frac{(-2x)^{\frac{11}{2}}}{64}+\cdots \\
\end{align*}
as one of the jets.  It can be shown that the the possible values of $c_4$ are equal to the roots of $\frac{(-2)^{\frac{11}{2}}}{64}$, and our solutions for $c_4$ seem to be imaginary.  However, $x^{\frac{11}{2}}$ forms a cusp, and the fact that there is a $-x$ in the jet simply means that the cusp opens left, so when $x$  is negative, the coefficient of the $\frac{11}{2}$ term in the jet is real.

\end{document}